\documentclass{article}

\usepackage{PRIMEarxiv}

\usepackage[noadjust]{cite}
\usepackage{multicol}

\usepackage[utf8]{inputenc}
\usepackage{eurosym}
\usepackage{bm}\usepackage{subfig}
\usepackage{amsmath,mathabx}
\usepackage{multirow}
\usepackage{booktabs}
\usepackage[capitalise]{cleveref}
\usepackage{float}
\usepackage[version=4]{mhchem}
\usepackage{siunitx}
\usepackage{caption,bbold}
\usepackage{longtable,tabularx}
\usepackage{algorithm}
\usepackage{algorithmic}
\usepackage{multirow}
\usepackage{url}
\usepackage{tikz}
\usepackage{mathtools}
\usepackage{cuted}
\usepackage{float}

\usepackage{float}	 											

\usepackage{setspace} 											
\usepackage{color}
\usepackage{caption}											
\newcommand\newl[2]{ 											
 	\expandafter\newlength\csname #1\endcsname
 	\expandafter\setlength\csname #1\endcsname{#2}}
\newcommand\setl[2]{ 											
 	\expandafter\setlength\csname #1\endcsname{#2}}

\newl{figh}{.0\columnwidth}										
\newl{figw}{.0\columnwidth}										
\newl{boxh}{.0\columnwidth}										
\newl{boxw}{.0\columnwidth}										
\newl{tcw}{.0\columnwidth}										

\newl{subfigskip}{.6em}											
\newl{plotvsep}{.8em}											

\usepackage{pgfplots}
\usepgfplotslibrary{groupplots}
\usetikzlibrary{positioning}

\pgfplotsset{
	compat=newest,
	plot coordinates/math parser=false,
}

\pgfdeclarelayer{annotation}
\pgfdeclarelayer{grid}
	\pgfsetlayers{grid,main,annotation}
\newcommand{\addimagegrid}{1}									
\newcommand{\imagegrid}{										
	\begin{pgfonlayer}{grid}
		\draw[black!5, ultra thin, step=1,] (img.south west) grid (img.north east);
	 	\draw[black!15, ultra thin, step=10,] (img.south west) grid (img.north east);
	 	\foreach \x in {0,10,...,100} { \node [below] at (\x,0) {\x};}
	 	\foreach \y in {0,10,...,100} { \node [left] at (0,\y) {\y};}
	\end{pgfonlayer}
}		

\usepackage[nopostdot,acronym,nogroupskip,nonumberlist]{glossaries}

\usepackage{xcolor}
\graphicspath{{Figures/}}

\usepackage{mathtools}

\makeglossaries
\newacronym{vlt}{VLT}{very low-thrust}
\newacronym{sa}{SA}{short arc}
\newacronym{ar}{AR}{Admissible Region}
\newacronym{np}{NP}{non-polynomial}
\newacronym{ads}{ADS}{Automatic Domain Splitting}
\newacronym{nlp}{NLP}{Nonlinear Programming Problem}
\newacronym{sqp}{SQP}{Sequential Quadratic Programming}
\newacronym{da}{DA}{Differential Algebra}
\newacronym{soc}{SOC}{Second Order Cone}
\newacronym{ocp}{OCP}{optimal control problem}
\newacronym{socp}{SOCP}{Second-Order Cone Programming}
\newacronym{dta}{DTA}{Defence Agency Technology}
\newacronym{rso}{RSO}{Resident Space Object}
\newacronym{rtn}{RTN}{Radial, Transverse, and Normal}
\newacronym{leo}{LEO}{low Earth orbit}
\newacronym{geo}{GEO}{geostationary Earth orbit}
\newacronym{tle}{TLE}{Two-Line Element}
\newacronym{uct}{UCT}{Uncorrelated Track}
\newacronym{dmc}{DMC}{dynamic model compensation}
\newacronym{snc}{SNC}{state noise compensation}
\newacronym{rms}{RMS}{root-mean-square}
\newacronym{cdm}{CDM}{control distance metric}
\newacronym{ewsk}{EWSK}{east-west station-keeping}
\newacronym{cut}{CUT}{Conjugate Unscented Transform}
\newacronym{lt}{LT}{low-thrust}
\usepackage[nopostdot,acronym,nogroupskip,nonumberlist]{glossaries}
\makeglossaries
\newacronym{ssa}{SSA}{space situational awareness}
\newacronym{iod}{IOD}{initial orbit determination}
\newacronym{pod}{POD}{precise orbit determination}
\newacronym{od}{OD}{orbit determination}
\newacronym{coe}{COE}{classical orbital elements}
\newacronym{cc}{CC}{cartesian coordinates}
\newacronym{mc}{MC}{Monte Carlo}
\newacronym{os}{OS}{Orbit Set}
\newacronym{vsa}{VSA}{very-short arc}
\newacronym{TFRM}{TFRM}{Telescope Fabra ROA Montsec}
\newacronym{asr}{ASR}{Admissible States Region}
\newacronym{par}{PAR}{probabilistic admissible region}
\newacronym{car}{CAR}{constrained admissible region}
\newacronym{vd}{VD}{virtual debris}
\newacronym{sd}{SD}{space debris}
\newacronym{otoa}{OTOA}{observation-to-observation association}
\newacronym{otta}{OTTA}{observation-to-track association}
\newacronym{ttta}{TTTA}{track-to-track association}
\newacronym{gmm}{GMM}{Gaussian mixture model}
\newacronym{ps}{PS}{point-wise sampling}
\newacronym{neo}{NEO}{near-Earth orbits}
\newacronym{tps}{TPS}{truncated power series}
\newacronym{ci}{CI}{confidence interval}
\newacronym{us}{US}{United States}
\newacronym{eu}{EU}{European Union}
\newacronym{ssn}{SSN}{space surveillance network}
\newacronym{sst}{SST}{space surveillance and tracking}
\newacronym{eusst}{EUSST}{European space surveillance and tracking}
\newacronym{JspOC}{JspOC}{Joint Space Operations Center}
\newacronym{USSTRATCOM}{USSTRATCOM}{U.S. Strategic Command}
\newacronym{fov}{FOV}{field of view}
\newacronym{esa}{ESA}{European Space Acency}
\newacronym{nasa}{NASA}{National Aeronautics and Space Administration}
\newacronym{ur}{UR}{Universidad de la Rioja}
\newacronym{vo}{VO}{Virtual Observatory}
\newacronym{meo}{MEO}{medium-Earth orbit}
\newacronym{heo}{HEO}{highly-elliptical orbit}
\newacronym{spice}{SPICE}{Spacecraft Planet Instrument C-matrix Events}
\newacronym{aida}{AIDA}{Accurate Integrator for Debris Analysis}
\newacronym{eci}{ECI}{Earth-centered inertial}
\newacronym{mssr}{MSSR}{Monopulse Surveillance Secondary Radar}
\newacronym{bssr}{BSSR}{Bistatic Surveillance Secondary Radar}
\newacronym{pie}{PIE}{parametric implicit equation}
\newacronym{ode}{ODE}{ordinary differential equation}
\newacronym{grv}{GRV}{Gaussian random variable}
\newacronym{as}{AS}{Admissible States}
\newacronym{gto}{GTO}{geostationary transfer orbit}
\newacronym{eom}{EoM}{Equations of Motion}
\newacronym{ut}{UT}{unscented transformation}
\newacronym{pc}{PC}{polynomial chaos}
\newacronym{stt}{STT}{state transition tensor}
\newacronym[\glslongpluralkey={state transition matrices}]{stm}{STM}{state transition matrix}
\newacronym{srp}{SRP}{solar radiation pressure}
\newacronym{pfe}{PFE}{Plank-Fokker equation}
\newacronym{daiod}{DAIOD}{Differential Algebra Initial Orbit Determination}
\newacronym{mee}{MEE}{modified equinoctial elements}
\newacronym{kep}{KEP}{Keplerian classical Elements}
\newacronym{pf}{PF}{Particle Filter}
\newacronym{ls}{LS}{Least Square}
\newacronym{pmd}{PMD}{Post Mission Disposal}
\newacronym{usssn}{US SSN}{US Space Surveillance Network}
\newacronym{aup}{AUP}{Association and Uncertainty Pruning}
\newacronym{mcmc}{MCMC}{Monte Carlo Markov Chain}
\newacronym{bfgs}{BFGS}{Broyden-Fletcher-Goldfarb-Shanno}
\newacronym{fom}{FoM}{figure of merit}
\newacronym{iet}{IET}{Ideal Elements (Time)}
\newacronym{ieT}{IE*}{Ideal Elements (\(\theta^{*}\))}
\newacronym{ie}{IE}{Ideal Elements}
\newacronym{ad}{AD}{Arbitrary Direction}
\newacronym{ss}{SS}{Subset Simulation}
\newacronym{tsa}{TSA}{too-short arc} 
\newacronym{aod}{AOD}{accurate orbit determination}
\newacronym{dav}{DAvector}{differential algebraic vector}
\newacronym{rmse}{RMSE}{root mean square error}
\newacronym{se}{SE}{standard estimate}
\newacronym{naif}{NAIF}{Navigation and Ancillary Information Facility}
\newacronym{et}{ET}{ephemeris time}
\newacronym{lov}{LOV}{Line of Variation}
\newacronym{ges}{GES}{gradient extremal surface}
\newacronym{iv}{IV}{initial value}
\newacronym{bv}{BV}{boundary value}
\newacronym{tpbvp}{TPBVP}{two-points boundary value problem}
\newacronym{utpbvp}{UTPBVP}{uncertain two-points boundary value problem}
\newacronym{tp}{TP}{true positive}
\newacronym{tn}{TN}{true negative}
\newacronym{fp}{FP}{false positive}
\newacronym{fn}{FN}{false negative}
\newacronym{mcs}{MCS}{Monte Carlo simulation}
\newacronym{ge}{GE}{Gradient Extremal}
\newacronym{mtt}{MTT}{multi-target tracking}
\newacronym{AIUB}{AIUB}{Astronomisches Institut - Universitat Bern}
\newacronym{GSTP}{GSTP}{General Support Technology Programme}
\newacronym{ZimSMART}{ZimSMART}{Zimmerwald SMall Aperture Robotic Telescope}
\newacronym{dace}{DACE}{Differential Algebra Computing Engine}
\newacronym{dals}{DALS}{Differential Algebra Least Squares}
\newacronym{dof}{DoF}{degrees of freedom}
\newacronym{lp}{LP}{Linear Programming}
\newacronym{sdp}{SDP}{Semidefinite programming}
\makeglossaries

\newcommand\doublerulefill{\leavevmode\leaders\vbox{\hrule width .1pt\kern1pt\hrule}\hfill\kern0pt }

\graphicspath{{media/}}     

\title{Detection and estimation of spacecraft maneuvers for catalog maintenance                                                                                                                                                                       
\thanks{© 2022 IEEE. Personal use of this material is permitted. Permission from IEEE must be obtained for all other uses, in any current or future media, including reprinting/republishing this material for advertising or promotional purposes, creating new collective works, for resale or redistribution to servers or lists, or reuse of any copyrighted component of this work in other works.} 
}

\author{
  Laura~Pirovano,  Roberto~Armellin \\
  Te P\=unaha \=Atea - Space Institute \\
  The University of Auckland \\
  Auckland, New Zealand\\
  \texttt{\{ laura.pirovano, roberto.armellin\}@auckland.ac.nz} 
}

\begin{document}
\maketitle

\begin{abstract}
Building and maintaining a catalog of resident space objects involves several tasks, ranging from observations to data analysis. Once acquired, the knowledge of a space object needs to be updated following a dedicated observing schedule. Dynamics mismodeling and unknown maneuvers can alter the catalog's accuracy, resulting in uncorrelated observations originating from the same object. Starting from two independent orbits, this work presents a novel approach to detect and estimate maneuvers of resident space objects, which allows for correlation recovery. The estimation is performed with successive convex optimization without a-priori assumption on the thrust arcs structure and thrust direction.
\end{abstract}


\begin{multicols}{2}
\printglossary[type=\acronymtype,title={List of Acronyms}]
\end{multicols}

\section{Introduction}
The space population is rapidly increasing. As of  June \num{2022}, more than \num{45,000} objects larger than \SI{10}{cm} are being tracked in space\footnote{Source: \url{www.space-track.org}}, and the active space population has doubled in the last four years. With the current development and deployment of very large constellations, it is clear that the safety and stability of operations need to be improved to allow long-term sustainability of space activities. To do so, an accurate catalog of all resident space objects is fundamental, which requires initialization of newly discovered tracks \cite{Pirovano2020b,Tao2020,Losacco2022}, correlation between tracks \cite{Pirovano2021}, and linkage of new observations to known orbits \cite{Pirovano2020a}. Once initialized, the catalog needs to be maintained through careful scheduling for object tracking. However, operational spacecraft perform frequent maneuvers, usually dictated by mission requirements or collision avoidance, potentially creating double entries in the catalog with post-maneuver measurements and thus hindering the data association problem. Recovering data association with the assumption of unknown maneuvers is thus fundamental to keep the catalog up-to-date and will be the focus of this manuscript. \\
\\

Past efforts in maneuver detection started from linked orbits to analyze the performance of satellites. Patera \cite{Patera2008} found maneuvers by identifying anomalous deviations between the published values of an orbital parameter and the corresponding values given by polynomial fitting of the moving window curve fit technique. Lemmens and Krag \cite{Lemmens2014} addressed maneuver detection for \glspl{leo} detecting space events based on analyzing their \glspl{tle} history. However, the most challenging aspect of maneuver detection is understanding whether two \glspl{uct} can be linked using a maneuver, thus regaining data association through the assumption of a maneuver, rather than detecting maneuvers between two orbits known to be the same object.\\
The most studied methods for dealing with uncertain dynamics in the state estimation process include the use of \gls{snc}, that is the addition of white noise to the process, and \gls{dmc},  which yields a deterministic acceleration term as well as a purely random term \cite{Tapley}.  However, they are not suitable for the problem at hand. Indeed, while effectively preventing divergence due to mismodeling, adding process noise does not provide a method for estimating the mismodeling or detecting its presence. On the other hand, \Gls{dmc} requires a significant amount of tuning, and appending dynamics parameters to the state requires a known model for those modeled dynamics \cite{Lubey2015}.
\\Lately, optimal \glspl{cdm} approaches \cite{Holzinger2012,Lubey2015,Singh2012} were shown to be suitable for detecting and characterizing maneuvers of \glspl{uct} in the absence of full orbit solutions. In \cite{Holzinger2012} it is assumed that the object which realizes the measurement with the least energy consumption, in linearized dynamics is the most likely originator. Furthermore, it gives a way to determine the probability that the nominal cost (a possible maneuver) is greater than the system uncertainty, hence considering non-perfect data association, with Gaussian assumption. In this way, the approach can reconstruct maneuvers with no a-priori knowledge. In \cite{Lubey2015} Lubey extended the \gls{cdm} approach to include several perturbations and adapted it to estimate unmodeled dynamics  - not only maneuvers. However, Singh et al. in \cite{Singh2012} showed that an energy-based metric rather than a propellant consumption-based one significantly deteriorates the quality of the reconstructed maneuvers. To avoid this issue, they formulated a \gls{nlp} problem, a step that affects their approach's efficiency and robustness. In summary, \glspl{cdm} approaches provide data-driven capabilities for detecting and characterizing unknown maneuvers. However, there is overall a lack of a methodology meeting all the following requirements: 
\begin{itemize}
    \item use \(\Delta\)V as a metric to enable a more reliable indicator of track correlation and a more realistic maneuver estimation,
    \item avoid linearized and/or simplified dynamics to potentially handle large uncertainties and data sparsity,
    \item work for arbitrary statistics,
    \item is numerically efficient and robust to be used in routine operations. 
\end{itemize}
This is the gap that this manuscript covers: following the definition of a reference trajectory with accurate dynamics, any deviation from it is propagated, and a minimum-fuel path is found through (successive) convex optimization. Focus is posed on the choice of coordinates for which a single iteration is enough to reach convergence. Within a single convexification, indeed, polynomial run time and global convergence are ensured. Two methods are then introduced to transfer the state uncertainty into maneuver uncertainty: the Mahalanobis distance approach and the \gls{cut}, which respectively deal with the detection and the estimation parts.\\
\\
The paper is organized as follows. \Cref{sec::dynman} describes the mathematical model for the fuel-optimal maneuver pattern to fit two states at different epochs, where the \gls{ocp} is transformed from a \gls{nlp} to a \gls{socp} one. \Cref{sec::cut}  introduces the Mahalanobis distance and \gls{cut} approaches,  defining the different constraints to follow the two methods. The full optimization algorithm is shown in \cref{sub:fullprob}. Results are shown in \Cref{sec::results}, where maneuvers are reconstructed for a synthetic maneuver, for a routine \gls{ewsk} maneuver in \gls{geo} and for a \gls{lt} orbit raising of a highly-elliptical orbit. Lastly, conclusions are drawn in \Cref{sec::conclusions}.

\section{Dynamics of a maneuvering object}
\label{sec::dynman}
This section explains the dynamical models and approaches followed to build the optimization problem. \Cref{sub::nonlin} introduces the formulation of a nonlinear \gls{ocp}, and defines the core constraints and objective of the maneuver detection problem. \Cref{sub::discnlp}  discretizes the problem, which is then convexified in \Cref{sub::convex}. 

\subsection{Formulation of the nonlinear \gls{ocp}}\label{sub::nonlin}
The problem of finding maneuver patterns that connect two defined states (\(x_0, x_1\)) given the time of flight (\(\Delta t\)) is a an \gls{ocp} with boundary constraints. Being the two available states the result of two independent \glspl{od}, neighboring values of the reference states can be explored within a set \(\mathcal{X}\), hence introducing uncertain boundary constraints. Assuming operators implement optimal maneuvers, the following nonlinear \gls{ocp} with boundary value conditions can be formulated, with the task of minimizing propellant use:
\begin{subequations}
\label{eq::fullNLPprob}
\begin{align}
   \min_u \quad & J( x, u, p \left|\right. x_0, x_1, t_0, t_f )\label{eq::J}\\
   \mbox{s.t.} \quad  & \dot{x}(t) = f(x, u, p,t), \label{eq::dyn}\\
  & \left(x(t),u(t), t \right) \in \mathcal{C}^i(t),& \forall t \in \mathbb{R}_{\Delta t}, i \in \mathbb{N}_0 \label{eq::cont}\\
  &   x_{t_0} \in \mathcal{X}_{t_0}, \quad  x_{t_f}\in \mathcal{X}_{t_f}.\label{eq::var}
\end{align}
\end{subequations}
The cost in \cref{eq::J} describes the goal of the mission to minimize the use of propellant, given the boundary conditions in terms of state and time. \Cref{eq::dyn} describes the system dynamics, which is composed of natural motion, spacecraft control \(u\) and any other parameters \(p\).  \Cref{eq::cont} ensures continuity or higher differentiability of the trajectory. Lastly, \cref{eq::var} bounds the initial and final states' variation. Each equation of Problem \ref{eq::fullNLPprob} is now described. 
\paragraph{Minimization function (\cref{eq::J})} the typical minimization function of a minimum fuel \gls{ocp} is
\begin{equation}
\label{eq::jnonlin}
\Delta V = \int_{t_0}^{t_f} \|\bm{u}(\tau)\| d\tau,
\end{equation}
that is the total thrusting \(\Delta V\) over the arc \(T=[t_0,t_f]\).
\paragraph{State dynamics (\cref{eq::dyn})} It is fundamental that the natural motion of the body is accurately described, to avoid overloading the control with unmodeled known dynamics, resulting in a non-realistic \(\Delta V\) profile. For this reason, the dynamics are described with the numerical propagator \gls{aida} \cite{Morselli2014}, which implements the following acceleration vector
\begin{equation}
    \bm{a}_{\textrm{AIDA}} = \bm{a}_\Earth + \bm{a}_\textrm{drag} + \bm{a}_\Sun + \bm{a}_\Moon  + \bm{a}_\textrm{SRP} = \bm{f} \left( \bm{x}, t\right),
\end{equation}
whose components are:
\begin{itemize}
\item[$\bm{a}_\Earth$]  Earth's gravity potential up to order \(15\),
\item[$\bm{a}_\textrm{drag}$] drag with NRLMSISE-00 atmospheric density model, 
\item[$\bm{a}_{\Sun,\Moon}$] third-body perturbing accelerations - Sun and Moon - with NASA's SPICE\footnote{\url{https://naif.jpl.nasa.gov/naif/toolkit.html}} toolkit \cite{ACTON199665}, 
\item[$\bm{a}_\textrm{SRP}$] \gls{srp} with dual-cone shadow model. 
\end{itemize} 
The acceleration is expressed in \gls{eci} coordinates and the dynamics is integrated with a Runge-Kutta 7-8 scheme. The control is then added to the natural motion of the body. 
\paragraph{Boundary conditions (\cref{eq::var})} 
The boundary conditions are the results of two independent \glspl{od}. Assuming they can be modeled as multivariate Gaussian random variables, they can be defined by their mean, \(\bm{\mu}_{t_0}\) and  \(\bm{\mu}_{t_f}\), and covariance \(\Sigma_{t_0}\) and \(\Sigma_{t_f}\), which are positive definite matrices. It is then possible to bound the variation from the mean values with a \(\chi\)-square quantile of level \(\alpha\) and \num{6} \gls{dof}:
\begin{equation}
   \frac{1}{2} ( \bm{x}_{i}-\bm{\mu}_{t_j})^T \Sigma_{t_j}    ( \bm{x}_{i}-\bm{\mu}_{t_j})\le q_{\chi^2} (\alpha, 6) = \mathcal{M}, \quad  j=\{0,f\}, \, i=\{0,N\}.
   \label{eq::Mahala}
\end{equation}

Despite the problem being now fully defined, in most cases, the solution of such an optimal control problem requires time consuming numerical procedures \cite{Malyuta2021}. This is because traditional methods to solve this type of problem - indirect and direct - either result in \glspl{tpbvp} highly sensitive to initial guesses, or are  \gls{np} hard, which means that the amount of computation required to solve the problem will not be limited by a bound determinable a priori. Unknowable computational time and lack of assured algorithm convergence are the kinds of obstacles that would preclude its applications in aerospace engineering problems that demand reliable and rapid solutions \cite{Liu2017}. Instead, the following sections will show how to solve Problem \ref{eq::fullNLPprob}, by firstly discretizing the \gls{nlp} problem and then transforming it into a (successive) convex one. For a single convex optimization problem convergence is guaranteed to a global minimum in polynomial time.

\subsection{Discretizing the dynamics of the nonlinear \gls{ocp}}
\label{sub::discnlp}
In this section, Problem \ref{eq::fullNLPprob} will be transformed from infinite-dimensional to finite-dimensional by discretizing the trajectory into a finite number of nodes $N+1$, assuming that the control can only assume a fixed value within each segment, between two consecutive nodes. The initial node \(\bm{x}_0\) is then defined at the earliest epoch \(t_0\), while the last node  \(\bm{x}_N\) at the second epoch \(t_f\).  The optimization variables are then $9(N+1)$, being the state \(\bm{x}_k\) and the control \(\bm{u}_k\)  at each node $k$, while the minimization function becomes the sum of the contributions of the control over the segments. Problem \ref{eq::fullNLPprob} then becomes: 
\begin{subequations}
\label{eq::fullprobdisc}
\begin{align}
   \min \quad & \sum_{i=0}^N \|\bm{u}_i\| \cdot \Delta t \label{eq::Jdisc}\\
  \mbox{s.t.} \quad  & {\bm{f}}_i  = \bm{a}_{\Earth,i} + \bm{a}_{\textrm{drag},i} + \bm{a}_{3B,i}  + \bm{a}_{\textrm{SRP},i} + \bm{u}_i & i \in \mathbb{N}_{[0,N]} \label{eq::dyndisc}\\
 &  \bm{x}_{{i+1}} = \bm{x}_{{i}}+  \int_{i}^{{i+1}} {\bm{f}}_{i}(\bm{x},\bm{u}_i,t) \,dt & i \in \mathbb{N}_{[0,N]-1}\label{eq::cont2} \\
  & \frac{1}{2}\| F_{t_j} (\bm{x}_{i}-\bm{\mu}_{t_j})\| \le  \mathcal{M}, \mbox{\hspace*{1cm}}F_{t_j} = \sqrt{\Sigma_{t_j}}, &\quad  j=\{0,f\}, \, i=\{0,N\}.\label{eq::var2}
\end{align}
\end{subequations}

\subsection{Problem convexification}
\label{sub::convex}
The discretized problem in \cref{sub::discnlp} is still a \gls{nlp} problem, with \gls{np} computational time and lack of assured algorithm convergence. However, by carefully choosing the type of constraints, coordinates and minimization function, it is possible to obtain a convex problem, for which the following properties hold:
\begin{itemize}
    \item if a feasible solution exists, it is the global optimum;
    \item the run time complexity is polynomial;
    \item there is no need for an initial guess.
\end{itemize}
These properties ensure global convergence and efficiency, which is why many aerospace applications make use of it, from low-thrust trajectory optimization \cite{Morelli2022}, to asteroid hopping transfers designs \cite{Liu2021},  soft landings \cite{SoftLanding2013} and autonomous collision avoidance \cite{Armellin2021}, to name a few. It is to be noted that there are three types of convex optimization: \gls{lp}, \gls{socp}, and \gls{sdp}, with the first one being the most optimized and the latter being the one with the most modeling power. In return, the first may be insufficient to model complex problems, while for the latter existing algorithms still do not scale well to the problem size, meaning the solution speed deteriorates rapidly as the problem size increases. For this reason, \glspl{socp} are attractive because they hold a good balance between the other two types of convex optimization. \\
This section introduces the mathematical operations to turn Problem \ref{eq::fullprobdisc} into a \gls{socp} problem:
\begin{subequations}
\label{eq::convexprob}
\begin{align}
   \min \quad & \bm{c}^T \bm{x} \label{eq::Jconvex}\\
   \mbox{s.t.} \quad  & A \bm{x} = \bm{b} \label{eq::dyndisc2}\\
  & \| F_i \bm{x} + \bm{d}_i \| \le \bm{p}_i^T \bm{x} + q_i, & i\in\mathbb{N}_{[1:L]}\label{eq::cone}
\end{align}
\end{subequations}
where the objective function is linear and constraints are linear and second-order cones. \Cref{eq::cone} is called a rotated quadratic cone \(\mathcal{Q}_r^{n+1}\), where \(n\) is the dimensionality of \(\bm{x}\). Whenever \(\bm{d}\) and \(\bm{p}\) are null, the cone is quadratic, \(\mathcal{Q}^{n+1}\). The two are related by the following change of variables and notation:
\begin{equation}
  \begin{cases}
  t = \bm{p}_i^T \bm{x} + q_i \\ \bm{s} = F_i \bm{x} +\bm{d}_i 
  \end{cases} \Leftrightarrow \left(t,\bm{s}\right) \in \mathcal{Q}^{m+1} 
\end{equation}
where \(F_i\bm{x} \in \mathbb{R}^m\). Though at first look Problem \ref{eq::convexprob} does not look too different from Problem \ref{eq::fullprobdisc}, it still needs to go through some classical techniques, called \emph{equivalent transformation},  \emph{relaxation}, and \emph{successive linearization}.
\paragraph{Equivalent transformation \label{par::eqtransf}}
This transformation is called exact convexification \cite{Low2014,Liu2017} or lossless convexification \cite{ACIKMESE2011341} technique since the solution of the transformed problem is also the solution of the original problem. It is widely used to convexify objective functions that cannot be discretized as linear functions. Indeed, despite \cref{eq::Jdisc} being the discretization of \cref{eq::jnonlin}, they both still include the norm of the control vector in the optimization function, which is nonlinear. By adding a slack variable, a nonlinear optimization function can be turned into a linear optimization function and a quadratic constraint:
\begin{equation}
\label{eq::eqtransf}
    \begin{array}{ccc}
    \cref{eq::Jdisc}&  \equiv & \begin{array}{ll}
        \min &\displaystyle\sum_i u_i \cdot \Delta t \\
        \mbox{s.t.} & \|\bm{u}_i\| = u_i
    \end{array}
    \end{array}
\end{equation}
 Once the control is discretized and transformed, it is possible to decide whether to  solve the optimization in acceleration, thus keeping \(\bm{u}_i\) as the optimization variable, or to add maneuvers as \(\Delta \bm{v}\)s at each node, hence solving the optimization in velocity. These two approaches bring different advantages: the former allows for an accurate representation of \gls{lt} maneuvers but cannot describe impulses, while the latter allows for an easy description of both types of propulsion, with the assumption that the maximum allowed impulse at each node is equivalent to the resulting variation with constant \gls{lt} acceleration over two consecutive nodes: \(\Delta v_{i,i+1} = a_{S/C} \Delta t \). Due to the ease of description of both types of propulsion, the velocity formulation was chosen for this work. At each node, a possible maneuver is thus hypothesized, which is bounded by \(\Delta v_{MAX}\), which is the maximum available  maneuver at each node.
\paragraph{Relaxation \label{par::relax}} 
\Cref{eq::eqtransf} linearized the optimization function introducing a quadratic constraint. This is, however, non-convex but can be turned into a second-order cone constraint by relaxing the equality:
\begin{equation}
\label{eq::relax}
   \|\bm{u}_i\| = u_i \; \xrightarrow{Relax.} \; \|\bm{u}_i\| \le u_i \Leftrightarrow (u_i, \bm{u}_i) \in \mathcal{Q}^{4}, \quad \forall i \in \mathbb{N}_{[0,N]}.
\end{equation}
\Cref{eq::var} on the other hand is already relaxed and the following holds:
\begin{equation}
\label{eq::boundary}
       \left(2\mathcal{M},\sqrt{\Sigma_{t_j}} ( \bm{x}_i - \bm{\mu}_{t_j}) \right) \in \mathcal{Q}^{7}, \quad  j=\{0,f\}, \, i=\{0,N\}.
\end{equation}
Though the relaxation is introducing a  larger feasible set than the original problem, \cref{eq::relax} will always reach equality at the optimal point, since the goal is to minimize \(\|\bm{u}\|\), meaning the relaxation is also a lossless convexification technique.
\paragraph{Successive linearization \label{par::linearization}}
The previous techniques introduced slack variables and relaxed constraints without modifying the solution of the original problem. Successive linearization, instead, does introduce a simplification, in that the original nonlinear problem is substituted with a first-order Taylor expansion around the initial ballistic trajectory, referred to as \(\tilde{\bm{x}}\). Remembering that at each node there are \num{9} optimization variables from \cref{sub::discnlp} - plus \num{1} slack variable through the equivalent transformation - it is possible to discretize the influence of maneuvers and deviations from the ballistic trajectory with:
\begin{equation}
\label{eq::contdisc}
\begin{array}{ll}
 {\bm{x}}_{i+1} \approx  \tilde{\bm{x}}_{i+1} + \left[ R_i | M_i\right] \left[\begin{array}{c}
        \hspace*{-0.2cm}  \bm{x}_{i} -   \tilde{\bm{x}}_{i} \hspace*{-0.2cm}\\
        \hspace*{-0.2cm} \Delta \bm{v}_{i} -   \Delta \tilde{\bm{v}}_{i} \hspace*{-0.2cm}
    \end{array}\right], 
    \end{array}
\end{equation}
where \(i \in \mathbb{N}_{[0,N]}\), the initial reference velocity variation $\Delta \tilde{\bm{v}}_{i}$ is null, and \(R_i\in\mathbb{R}^{6x6}\), \(M_i\in\mathbb{R}^{6x3}\) are the \glspl{stm}, which respectively map deviations in the states and velocity to the next node. Due to linearization, iterations \(k\) may be  necessary to accurately solve the original problem: the optimization can thus be repeated following a standard successive convex optimization approach, where the new trajectory is convexified, until the solution of the optimization and the accurate forward propagation of the maneuvering orbit match to a prescribed accuracy, that is until two consecutive successive iterations yield the same optimization vector. However, successive convexification does not hold the same mathematical properties as a single convexification. Hence, it is also important to identify the state coordinates for which linearization holds for larger variations, such that the iterative process is not triggered.

\section{Handling uncertain constraints}
\label{sec::cut}
The optimization finds a deterministic path for this \gls{ocp}. However, the boundary constraints are uncertain, since the initial and final deviations to be found are not purely geometrical but hold a statistical meaning: the further from the mean states the final solution is, the less probable the path found is to happen. By simply constraining the initial and final deviations to be less than a value as in \cref{eq::var}, the optimizer will see the deviations as ``free propellant'' hence taking full advantage of them, but actually creating the least probable path, among the possible ones. At the optimum, the relaxed bound will, indeed, be exact.

\begin{figure}[ht]
\begin{center}
 \subfloat[Running optimizations with increasing value for the bound on state deviations.\label{fig::man_opt}]{\includegraphics[width= 0.45\columnwidth]{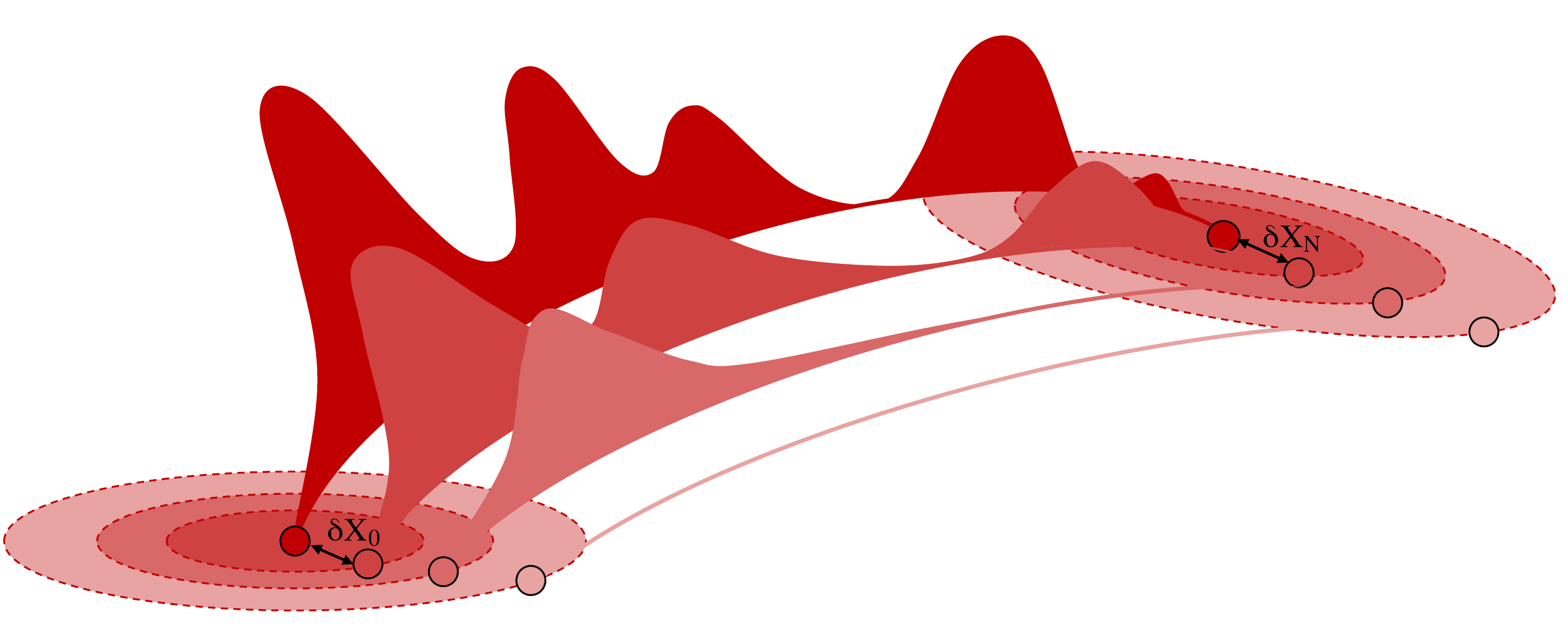}}
 \hfill
\subfloat[Running optimizations with predefined value of states deviations from \gls{cut}.\label{fig::man_CUT}]{\includegraphics[width= 0.45\columnwidth]{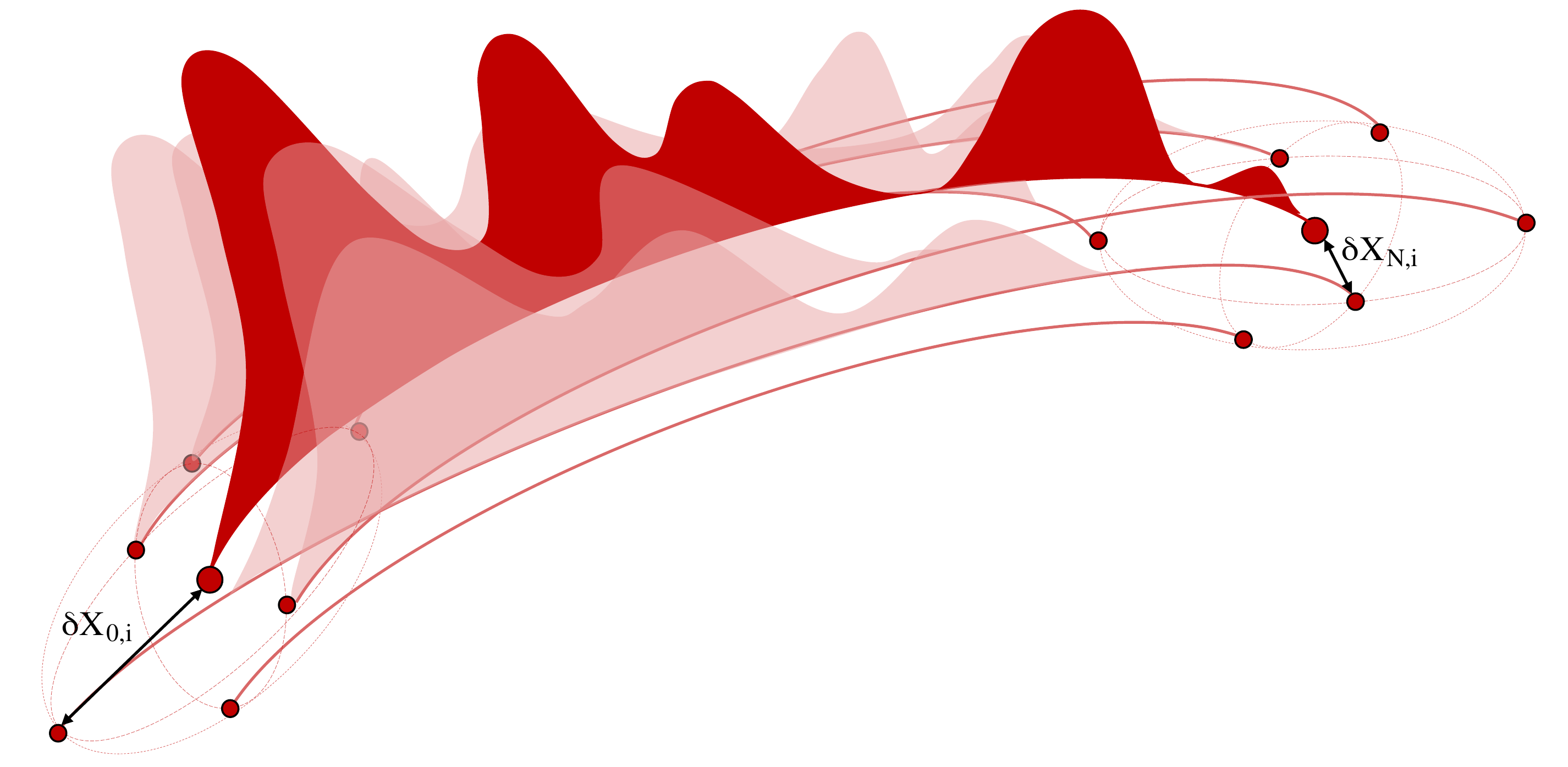}}
\caption{Dealing with uncertain constraints inside (left) and outside (right) the optimization routine, by respectively successively constraining the Mahalanobis distance value and using the \gls{cut}-4 algorithm.}
\end{center}
\end{figure}

Two different approaches can then be followed, depending on the information sought. If an indication on whether the maneuver has happened or not is wanted, it is possible to set a bound on the maximum allowed Mahalanobis distance. The optimization will then return the minimum \(\Delta v\) required to link the states with the chosen confidence. With just a few runs with different bounds it is then possible to obtain an indication on the probability that a maneuver has happened. This is sketched in \cref{fig::man_opt}. However, the optimization always finds the optimum on the ellipsoid, hence hindering the a-posteriori statistical analysis. For this reason, if an accurate statistics of the maneuver is sought, the initial and final states variations cannot be determined within the optimization but need to be defined beforehand,  shown in \cref{fig::man_CUT}. The fourth-order \gls{cut} \cite{Adurthi2012} was chosen to do so: by running the optimization \(N^2 + 2N +1\) times, where \(N=12\) is the dimension of the problem, it is possible to reconstruct a-posteriori the first four moments - mean, variance, skewness, and kurtosis - of the final maneuver profile distribution, thus allowing for the analysis of the effect of uncertainty in the states on the overall maneuver existence and estimation. Sigma points drawn from the initial and final covariance determine the initial and final deviations, so that the cone constraint in \cref{eq::boundary} is substituted by:
\begin{equation}
\label{eqCUT}
    \left[\begin{array}{c}
        \hspace*{-0.2cm}\bm{x}_0 - \bm{\mu}_{t_0} \hspace*{-0.2cm}\\
        \hspace*{-0.2cm} \bm{x}_N - \bm{\mu}_{t_f} \hspace*{-0.2cm}
    \end{array}\right]_c = \Delta \textbf{X}^{12x1}_{CUT4,c} \quad c\in \mathbb{N}_{[1, N^2 + 2N + 1]}
\end{equation}
where for each \(c\) the optimization is run. The superiority in efficiency of the \gls{cut} with respect to a simple \gls{mc} run and other cubature methods is detailed in \cite{Adurthi2018}.  Note that also with this construction  the ballistic trajectory is only computed once and holds for all \gls{cut} iterations, since deviations are linearly propagated through \glspl{stm}.

\section{Finalization of the \gls{socp} problem}
\label{sub:fullprob}
\begin{figure*}[!ht]
\centering
	\begin{tikzpicture}
	\node[inner sep=0pt, anchor=south west] (img) at (0,0) {\includegraphics[height=0.9\textwidth,angle=-90]{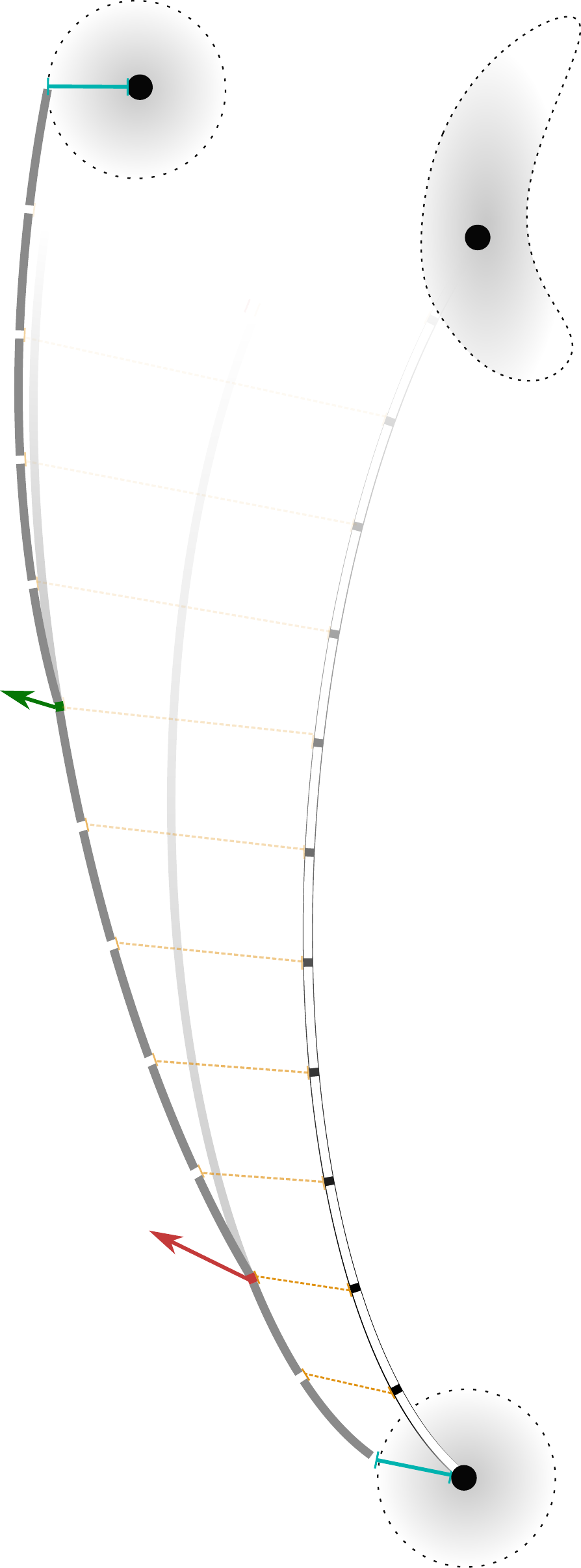}};
	\begin{scope}[ x={($0.01*(img.south east)$)}, y={($0.01*(img.north west)$)}, ]
		\ifnum\addimagegrid=0\imagegrid\fi
		\begin{pgfonlayer}{annotation}

  		 	\node[inner sep=0pt,anchor=north] at (6,15) {\scriptsize $\bm{\mu}_{t_0},\Sigma_{t_0}$};
  		 	\node[inner sep=0pt,anchor=north] at (86,15) {\scriptsize prop($\bm{\mu}_{t_0},\Sigma_{t_0}$)};
  		 	\node[inner sep=0pt,anchor=north] at (94.5,73.2) {\scriptsize $\bm{\mu}_{t_f},\Sigma_{t_f}$};
 			\node[inner sep=0pt,anchor=south east] at (5.5,24) {\scriptsize $(\bm{x}_0- \bm{\mu}_{t_0})$};
 			\node[inner sep=0pt,anchor=south east] at (93.5,80) {\scriptsize $(\bm{x}_N- \bm{\mu}_{t_f})$};
 			\node[inner sep=0pt,anchor=south east] at (19,66) {\scriptsize $\Delta \bm{v}_i$};
 			\node[inner sep=0pt,anchor=south east] at (54,93) {\scriptsize $\Delta \bm{v}_j$};
		\end{pgfonlayer}
	\end{scope}	
	\end{tikzpicture}
\caption{Initial sampled trajectory (white line with black squares), initial deviation (blue) and mid-course maneuvers (red and green) to match final deviation (blue).}
\label{fig::traj}
\end{figure*}
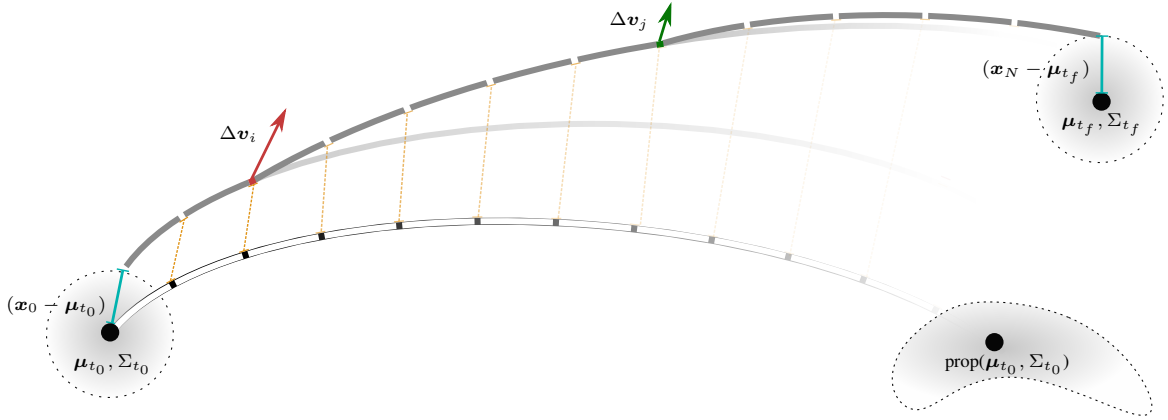

Now that the original problem is turned into a \gls{socp} and the uncertain nature of the bounds is taken care of, it is possible to enunciate the formulation to be solved. In order to make the description clearer, refer to \cref{fig::traj} for the remainder of the Section. Suppose two states - the outcome of an \gls{od}  \cite{LoSacco2021} - are available at two epochs, simplified with a black dot and a shaded ellipse, \(\left(\bm{\mu}_{t_0}, \Sigma_{t_0}\right)\), \(\left(\bm{\mu}_{t_f}, \Sigma_{t_f}\right)\). The state at the earliest epoch is propagated forward to the second one, determining a reference trajectory, white solid line. The propagated and determined states will not correlate with typical data association techniques \cite{Hill2008, Siminski2014, Pirovano2020b}, but correlation may be recovered assuming a maneuver has happened. The time interval is then discretized in N\(+1\) nodes, black squares on white line, and an impulsive maneuver \(\Delta \bm{v}_i\) is hypothesized at every node.
From \cref{eq::contdisc}, \glspl{stm} \(M_i\) and \(R_i\) are used to map the maneuvers' and state deviation effects forward in time. The cumulative effect of all maneuvers and initial deviation allows the state at \(\bm{x}_N\) to reach the desired location. The complete successive convex optimization problem with the \gls{cut}-4 approach is summarized in \cref{eq::fullprob2}. This problem can be solved by primal-dual interior-point methods. We use MOSEK through its MATLAB interface\footnote{\url{https://docs.mosek.com/9.0/toolbox/index.html}}, while \Gls{da}\footnote{\url{https://github.com/dacelib/dace}} is used to automatically compute the \glspl{stm}. Though no bounds on optimization variables are needed for the problem to converge, they are included to obtain realistic solutions. The \(\Delta v\)s are bounded for physical purposes, as stated before, so that the achievable maneuver can be obtained by the on-board propellant:
\begin{equation}
    \label{eq:deltavbound}
      0\le \Delta v_{i} \le \Delta v_{MAX}, \,|\Delta v_{il}| \le \Delta v_{MAX}\quad   \forall i \in \mathbb{N}_{[0,N]},\, \forall l = x,y,z,
\end{equation}
States, on the other hand, are bounded due to the linearization: a large state variation may introduce a significant error on the linearized dynamics, hence needing more iterations to converge:
\begin{equation}
    \label{eq:boundstate}
    \Delta x_{ij,MIN} \le x_{ij} - \tilde{x}_{ij} \le \Delta x_{ij,MAX},  \quad \forall i \in \mathbb{N}_{[0,N]}, \, \forall j \in \mathbb{N}_{[1,6]}
\end{equation}
It is expected that highly varying coordinates such as \gls{cc} take longer to converge to a solution than, for example, \gls{mee}. Indeed, examples in \cref{sec::results} will show that \gls{mee} are the preferred choice due to having only one fast variable and no singularities.

\begin{table*}[!ht]
\caption{Pseudo-code of the successive convex optimization algorithm for maneuver detection and estimation. The \emph{for}-loop is not performed when the Mahalanobis distance approach is taken.}
\label{eq::fullprob2}
\doublerulefill
\begin{equation*}
    \begin{array}{rrrlr}
       \multicolumn{4}{l}{\mbox{Initial guess} \quad\bm{x}_i^0 = \tilde{\bm{X}}_i^0, \quad \Delta \bm{v}^0_i = \Delta \tilde{\bm{v}}^0_i = \bm{0} \quad \forall i \in \mathbb{N}_{[0,N]}} & \emph{Ballistic trajectory}  \\[0.15cm]
    \multicolumn{4}{l}{ \mbox{for }c\in \mathbb{N}_{[1, N^2 + 2N + 1]} }\hfill{ } & \emph{Loop on CUT-4 sigma points}  \\[0.15cm]
    & \mbox{do} \hfill{ }\\[0.15cm]
  &  & \min & \sum_{i=0}^N \Delta v_i^k,\\[0.15cm]
& &   \mbox{s.t.} &  \cref{eq::contdisc}&   \emph{Continuity of trajectory}  \\[0.15cm]
&&& \cref{{eq::var2}} \mbox{ [No \emph{for}-loop],  or } \cref{eqCUT} & \emph{Initial and final states variation} \\[0.15cm]
&&& \cref{eq::relax} & \emph{Cone constraint on \(\Delta v\) norm}\\[0.15cm]
&&&  \cref{eq:deltavbound} & \emph{Bound on \(\Delta v\)}\\[0.15cm]
&&& \cref{eq:boundstate}& \emph{Bound on state variation} \\[0.15cm]
 &  \multicolumn{3}{l}{\mbox{while }\left(\|\bm{sol}^k - \bm{sol}^{k-1}\| > \epsilon\right)} & \emph{Successive convexification}\\
\mbox{end}
    \end{array}
\end{equation*}
\doublerulefill
\end{table*}

\section{Maneuver detection and estimation results}
\label{sec::results}
This section shows some test cases in which the exact reconstruction, the detection and lastly the estimation of maneuvers are performed. Firstly, the estimation of a synthetic maneuver is carried out, to show the \gls{lt} and impulsive reconstructions for perfect states. Then the Mahalanobis distance approach is shown for two cases taken from EUMETSAT publicly available dataset for an \gls{ewsk} maneuver and an non-maneuvering object. This will highlight the maneuver detection feature. Lastly, the \gls{cut} approach is used for the same \gls{ewsk} maneuver and for the \gls{lt} orbit raise from \gls{gto} to \gls{geo}, to showcase the maneuver estimation feature. The very large deviation from the reference trajectory will highlight the different behaviors of the coordinate systems. These very different orbits and maneuver types will showcase the versatility and tailorability of the approach. 

\subsection{Synthetic maneuver estimation}
\label{sec:sim}
\begin{figure}[ht]
\begin{center}
\includegraphics[width=0.5\columnwidth]{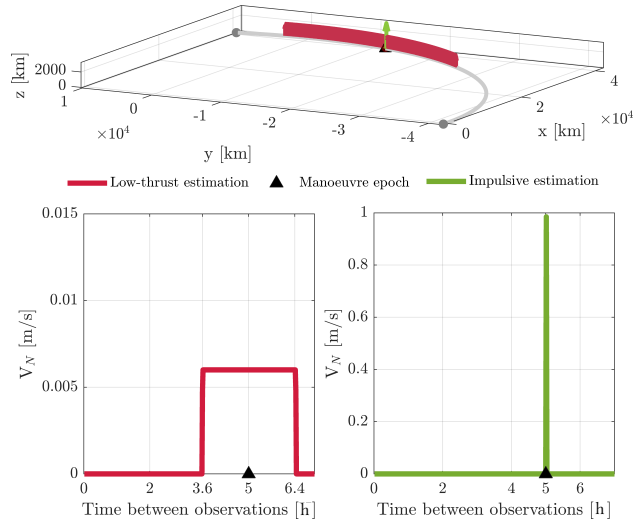}
\caption{Manoeuvre profile estimation of a simulated manoeuvre with \gls{lt} (left) and impulsive (right) assumption.}
\label{fig:sim}
\end{center}
\end{figure}

Maneuvers were initially simulated to test the algorithm, given known dynamics and perfect states. This was to test the accurate reconstruction of the maneuvers. \Cref{fig:sim} shows the reconstruction of a \num{1}~m/s out-of-plane maneuver performed $5$~h after the availability of the first state in keplerian dynamics. As shown in \cref{sec::dynman}, the initial trajectory dynamics is a plug-in, meaning any type of dynamics can be accommodated. In this case, keplerian dynamics was used to estimate the maneuver, so that it matched the dynamics of the object. The initial \gls{geo} trajectory was sampled every minute. Two reconstructions were performed to show the effect of the maximum allowed \(\Delta v\) at each node on the final profile: the unconstrained reconstruction (green) matches the impulsive input data, while the \gls{lt} reconstruction symmetrically spans around the maneuver epoch with a \(\Delta v_{MAX}= 6\)~{mm/s}. To obtain the desired maneuver, then, \num{167} small impulses were needed, summing up to \num{2.8}~h, which is the exact span estimated. 
\subsection{Mahalanobis distance approach}\label{res::mahala}
\begin{figure}[ht]
\subfloat[Maneuver detection run on two states a week apart with no maneuver involved. There exists a ballistic solution starting from \(68\%\) confidence. \label{fig::mahalaleft}]{\includegraphics[width=0.48\linewidth]{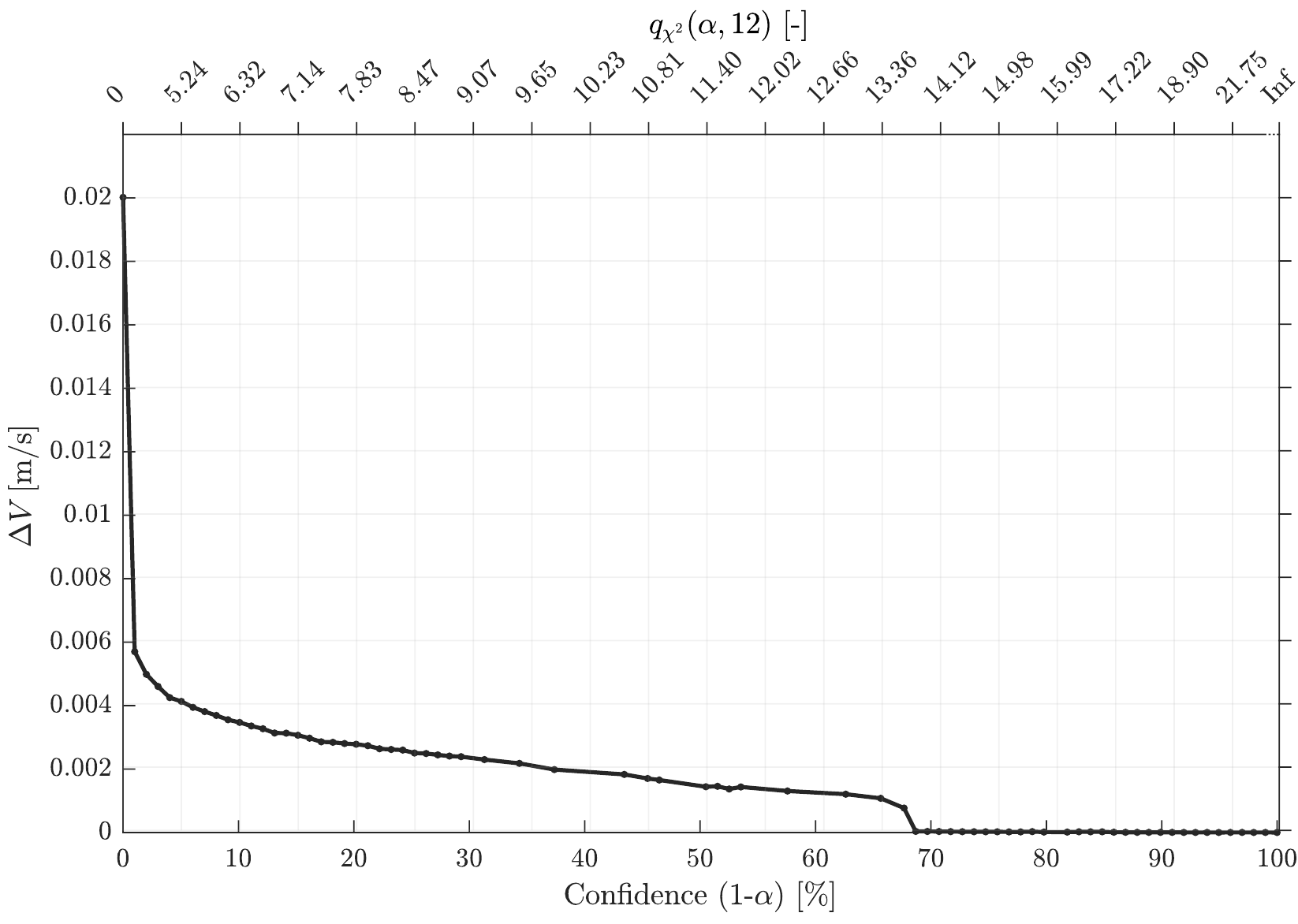}} 
\hfill
\subfloat[Maneuver detection run on two states a week apart with an \gls{ewsk} maneuver involved. A ballistic solution is never .\label{fig::mahalaright}]{\includegraphics[width=0.48\linewidth]{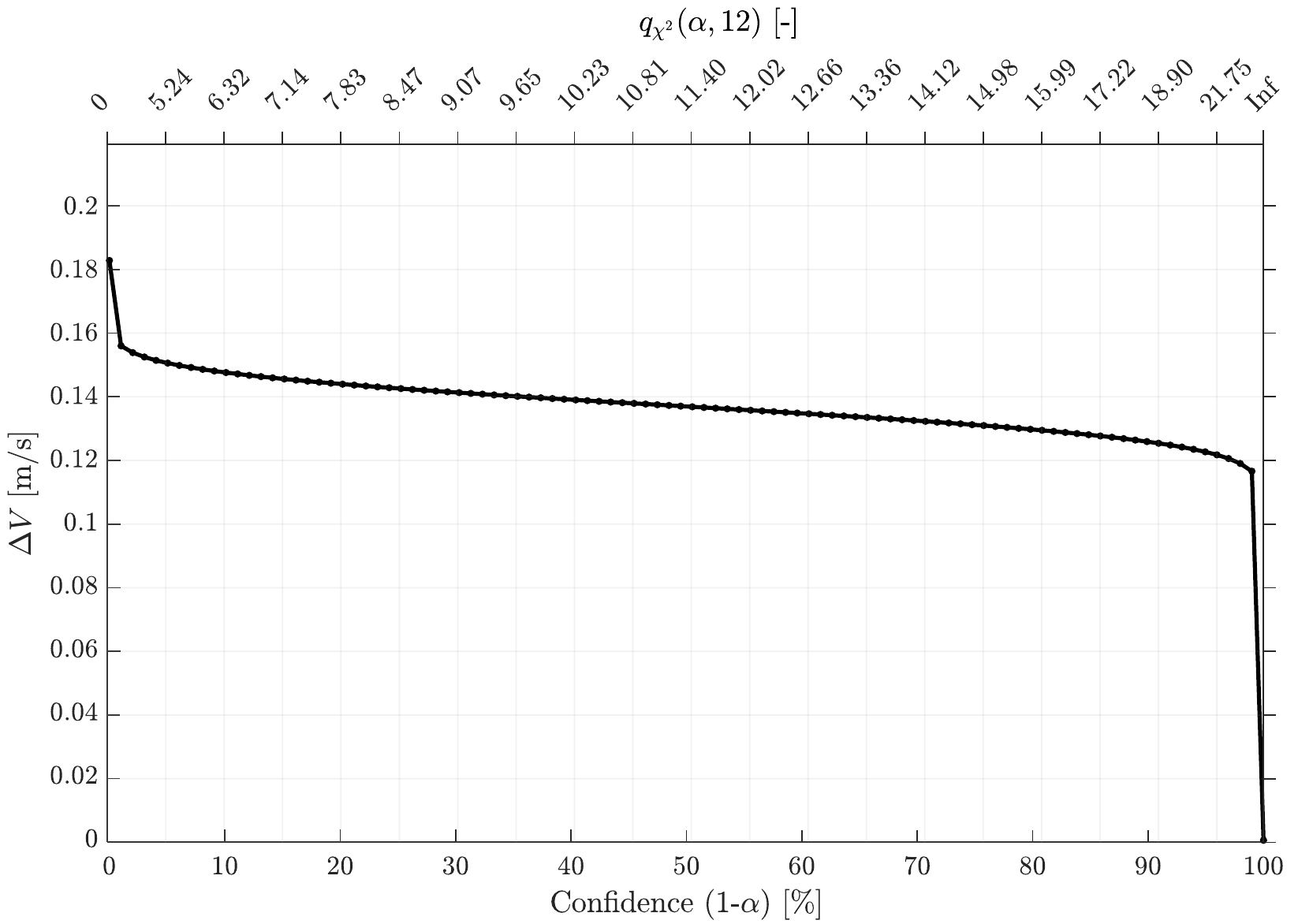}}
\caption{Evolution of minimum \(\Delta v\) required to connect two states depending on confidence level chosen.}
\label{fig::mahala}
\end{figure}

With the Mahalanobis distance approach, the states uncertanties are bounded by the ellipsoid defined in \cref{eq::Mahala}, where the value of \(\mathcal{M}\) depends on the confidence level \(\alpha\). The resulting maneuver profile represents the smallest maneuver achievable with the given confidence on the states. The larger the confidence, the larger is the search space around the mean state, hence the smaller the connecting maneuver can be. At the extreme, with \(\alpha = 0\) - hence \(100\%\) confidence - the associated \(\chi^2\) quantile is infinite, hence always finding a ballistic solution to connect the two orbits. \Cref{fig::mahala} shows the total maneuver magnitude depending on the confidence chosen for two different cases taken from EUMETSAT public dataset, one where an \gls{ewsk} maneuver is performed, and one where no maneuver is performed. Usually confidences of \(95\%\) or more are used in statistics to perform these kinds of analyses, hence showing that for \cref{fig::mahalaleft} there exists a ballistic solution which can connect the two orbits, while for \cref{fig::mahalaright} a residual maneuver of at least \(\Delta V = 0.12 \) m/s is needed to connect the two orbits with the same confidence. For comparison purposes, the accumulated \(\Delta v\) over a week for the \gls{srp} perturbation is \(\Delta V_{SRP} = 0.007 \) m/s. This is a very efficient way to detect maneuvers, as the \gls{stm} creation for \(50\) nodes takes 0.7 s and the optimisation takes \(25\) ms for the case at hand.

\subsection{\gls{cut} approach}
Differently from the Mahalanobis distance approach, the \gls{cut} approach can not only detect, but also estimate the maneuver profile and its higher order moments by running the optimization \(N^2 + 2N +1\) times. This means that the maneuver estimation comes at the expense of more computational power needed. \Cref{res::ewsk} shows the estimation of the same \gls{ewsk} maneuver analysed in \cref{res::mahala}, applying a further filter to the maneuver estimation process, that is the possibility to shrink the search window within a specific revolution. \Cref{res::mev2} instead estimates a very large \gls{lt} maneuver, hence introducing a large deviation from the reference trajectory. In this scenario, the choice of coordinates will greatly impact the computation of the maneuver. 
\subsubsection{\gls{ewsk} maneuver in \gls{geo}: reducing the maneuver window} \label{res::ewsk}
\begin{figure}[!ht]
\begin{center}
\includegraphics[width=0.7\columnwidth]{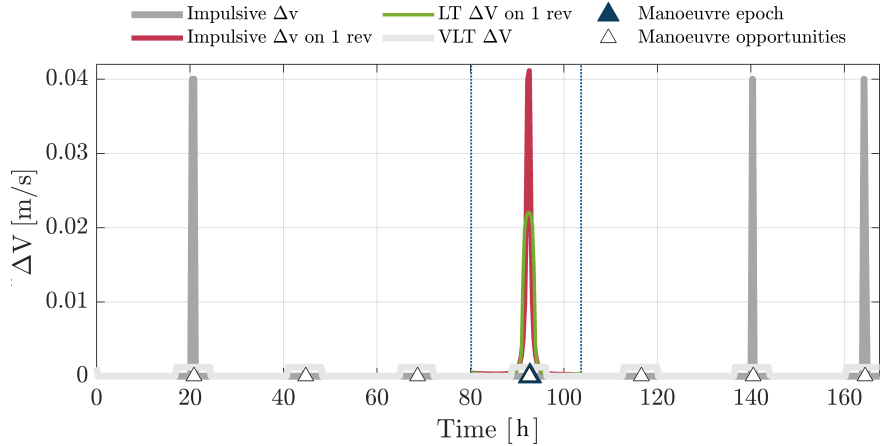}
\caption{\gls{ewsk} maneuver estimation on a seven days window and on a restricted one day window with impulsive and \gls{lt} assumption. Seven maneuver occasions are found, among which the real maneuver epoch. Estimation is performed on the assumption of perfect states.}
\label{fig:window}
\end{center}
\end{figure}

\noindent  EUMETSAT releases weekly updates of their METEOSAT satellites' states and covariances, including maneuvers types and times, if available. By estimating the maneuver, one can thus check whether the available data matches the reconstructed ones. However, having data for \gls{geo} satellites available every seven days, the same routine maneuver opportunity presents itself at each period, meaning a simple impulse may be split over seven different occasions, as shown in \cref{fig:window}, where the \gls{vlt} reconstruction greatly highlights this pattern. To exploit this information about the maneuver and simultaneously reduce the computational effort, the nodes were not evenly placed on the entire seven days window, but only defined on a portion of it. To do so, an estimate of the maneuver epoch was found by looking for the minimum distance between the forward propagation of \(\bm{\mu}_{t_0}\) and the backward propagation of \(\bm{\mu}_{t_f}\), then a symmetric interval was chosen around that epoch. The dashed column of \cref{fig:window} shows the new restricted window with the impulsive (red) and constrained (green) maneuver reconstruction.

\begin{figure}[!ht]
\begin{center}
\includegraphics[height=0.38\textheight]{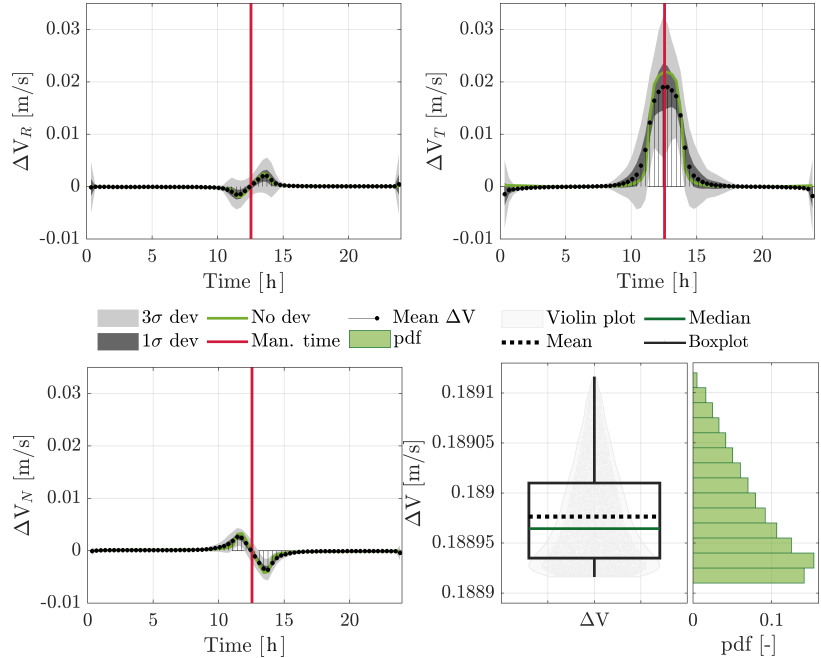}
\caption{\gls{ewsk} maneuver estimation for radial (R), tangential (T), and normal (N) components on a restricted one period window, with a-posteriori probability density reconstruction through the use of the \gls{cut}-4 algorithm.}
\label{fig:EWSK}
\vspace*{0.3cm}
\end{center}
\end{figure}

\noindent For the constrained reconstruction, \(\Delta v_{MAX} = 0.021\)~m/s was chosen for each node, due to lack of accurate data for the satellite.  It can be noted that the estimated peak matches with the real maneuver epoch.  While \cref{fig:window} was reconstructed with the assumption of perfect states, \cref{fig:EWSK} shows the complete reconstruction with a restricted window after the application of the \gls{cut}-4 algorithm. The green profile matches \cref{fig:window} solution being the first sigma point considered in the \gls{cut}-4. The total estimated maneuver was \(\Delta V= 0.18896\)~m/s, which fits in the typical range for \gls{ewsk} maneuvers, which sits between \num{0.05} and \num{0.2}~m/s \cite{larson1992space}. The bottom right plot in \Cref{fig:EWSK} shows the probability distribution, boxplot, and violin plot of the reconstructed maneuver magnitude. This was performed by fitting a Pearson distribution in the four moments estimated through the \gls{cut}-4. As can be seen, estimating the skewness \(s=0.7\) and kurtosis \(k=2.52\) parameters is fundamental to reconstruct the non-gaussian profile. This is especially evident for very small maneuvers where the left tail of the distribution is much smaller than the right tail. The cumulative distribution function of the estimated maneuver can then be used if competing reconstructions are estimated to find the most probable, as done in \cite{Holzinger2012}. The initial trajectory with \glspl{stm} creation took \(0.7\)~s (same process as for the Mahalanobis distance approach), each optimization took \(5-30\)~ms, for a total of \(22\)~s with parallelization in Matlab. The \gls{cut} approach thus takes advantage of the \gls{stm} creation and uses the output for all optimizations to cut the overall computational time. 

\subsubsection{\gls{lt} orbit raise maneuver from \gls{gto} to \gls{geo}: choosing the right coordinates} \label{res::mev2}
\begin{figure}[!ht]
\begin{center}
 \includegraphics[height=0.38\textheight]{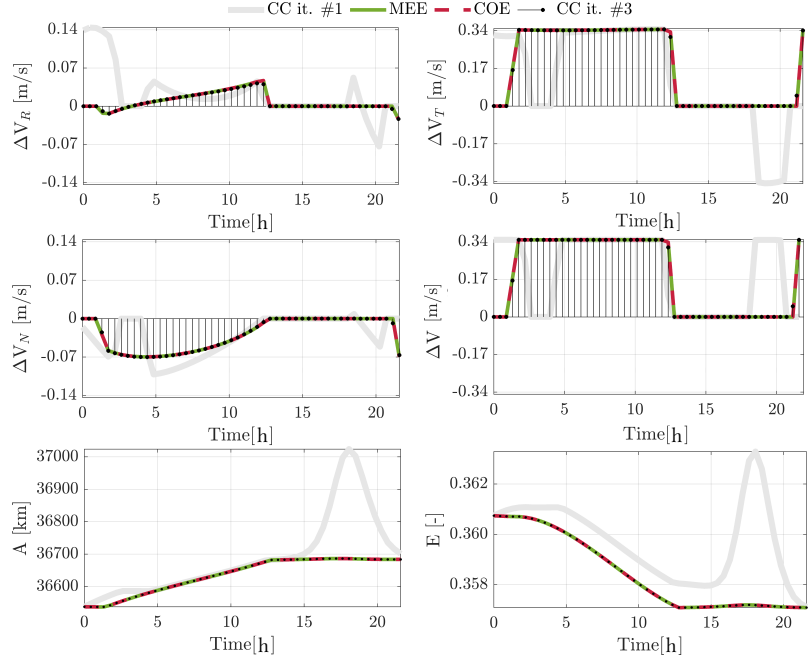}
\caption{\gls{lt} maneuver reconstruction with different coordinates. \gls{cc} need three successive convexifications to reach convergence, while \gls{coe} and \gls{mee} do not need any iterations. Non-convergence on first iteration is visible in the semimajor axis and eccentricity evolution (bottom). \label{fig::diffcoord}}
\end{center}
\end{figure}

For this test case we obtained tracking data of MEV-2 from \gls{dta}.States and covariances were then obtained through our in-house \gls{iod} algorithm \cite{Pirovano2020b}, the solution of which was then fed to our \gls{dals} algorithm \cite{LoSacco2021}. Low-thrust maneuvers were then estimated with a maximum acceleration of \(a=0.22\)~mm/s\(^2\). This case, differently from the \gls{ewsk} maneuver, involved a large maneuver on a highly elliptical orbit. It was thus of interest to understand the convergence of the algorithm depending on the coordinates used. For this reason, the maneuver was firstly reconstructed on the assumption of a perfect state for \gls{cc}, \gls{coe} and \gls{mee}. Results are shown in \cref{fig::diffcoord}. 
\begin{figure*}[!ht]
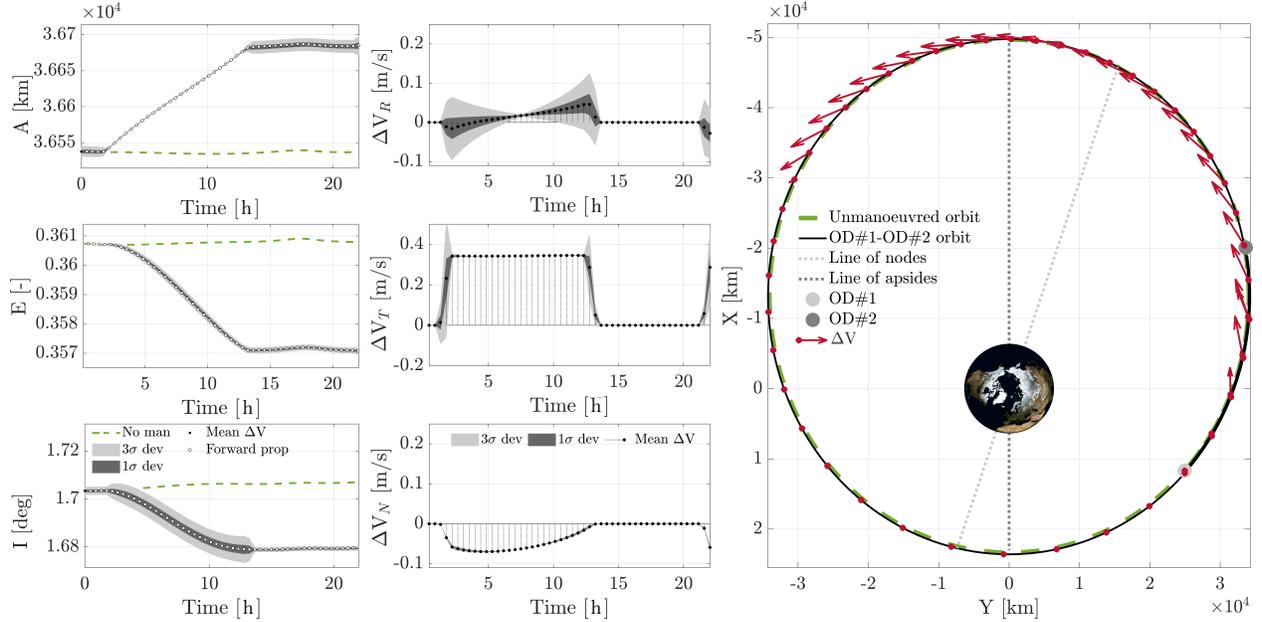

\begin{center}
\includegraphics[height=8.2cm]{DV12_50nodes_DTA.pdf}\hfill
\includegraphics[height=7.9cm]{DV_12_stem.pdf}\hfill
\includegraphics[height=8.2cm]{PQW_50nodes_DTA.pdf}
\caption{\gls{lt} maneuver estimation for an orbit raise by MEV-2 from \gls{gto} to \gls{geo}. The orbital elements variation (left) and maneuver profile in RTN components (middle) are reconstructed showing mean and \(3\sigma\) deviation through the \gls{cut}-4 algorithm. The maneuvers location along the orbit is shown on the right. Forward propagation of the obtained maneuver profile is also shown in the left plot.}
\label{fig::LTGTOGEO}
\end{center}
\end{figure*}
It can be seen that all coordinates reached convergence on the same maneuver pattern. However, \gls{coe} and \gls{mee} did not need any successive convexification, while \gls{cc} needed three iterations to reach convergence. The reconstructed patterns on the semimajor axis and eccentricity highlight the initial difficulty, where the linearization was not accurate.  This was expected, due to the fast varying nature of all components of \gls{cc}, as anticipated in \cref{sub::convex}. Despite all coordinates types reaching the same maneuver pattern, \gls{cc} performed the worst in terms of computation time, needing three successive convex optimizations to reach the same conclusion. \Cref{fig::LTGTOGEO} shows the complete reconstruction with elements evolution, maneuver profile, and orbital positioning of the maneuver obtained with \gls{mee}. Given that a change in the orbit shape and inclination happened, maneuvers were estimated around the line of apsides and the line of nodes, respectively, as expected. In the left plot, the validation is also shown: the white dots are the forward propagation implemented with the maneuver pattern estimated in the middle plot showing that the solution of the convexified problem also solved the original problem. The total maneuver was estimated to be \(\Delta\)V\(=9.1\)~m/s and took \(8.5\)~min: \(2.6\)~s for the trajectory and \glspl{stm} generation with \(50\) nodes, around \(0.37-0.83\)~s for each optimization , parallelized in Matlab for \(4121\) samples. The computational burden is thus due to the wish to estimate the \(\Delta V\) statistical properties. 

\begin{figure}[ht]
  \begin{center}
    \includegraphics[width=0.35\columnwidth]{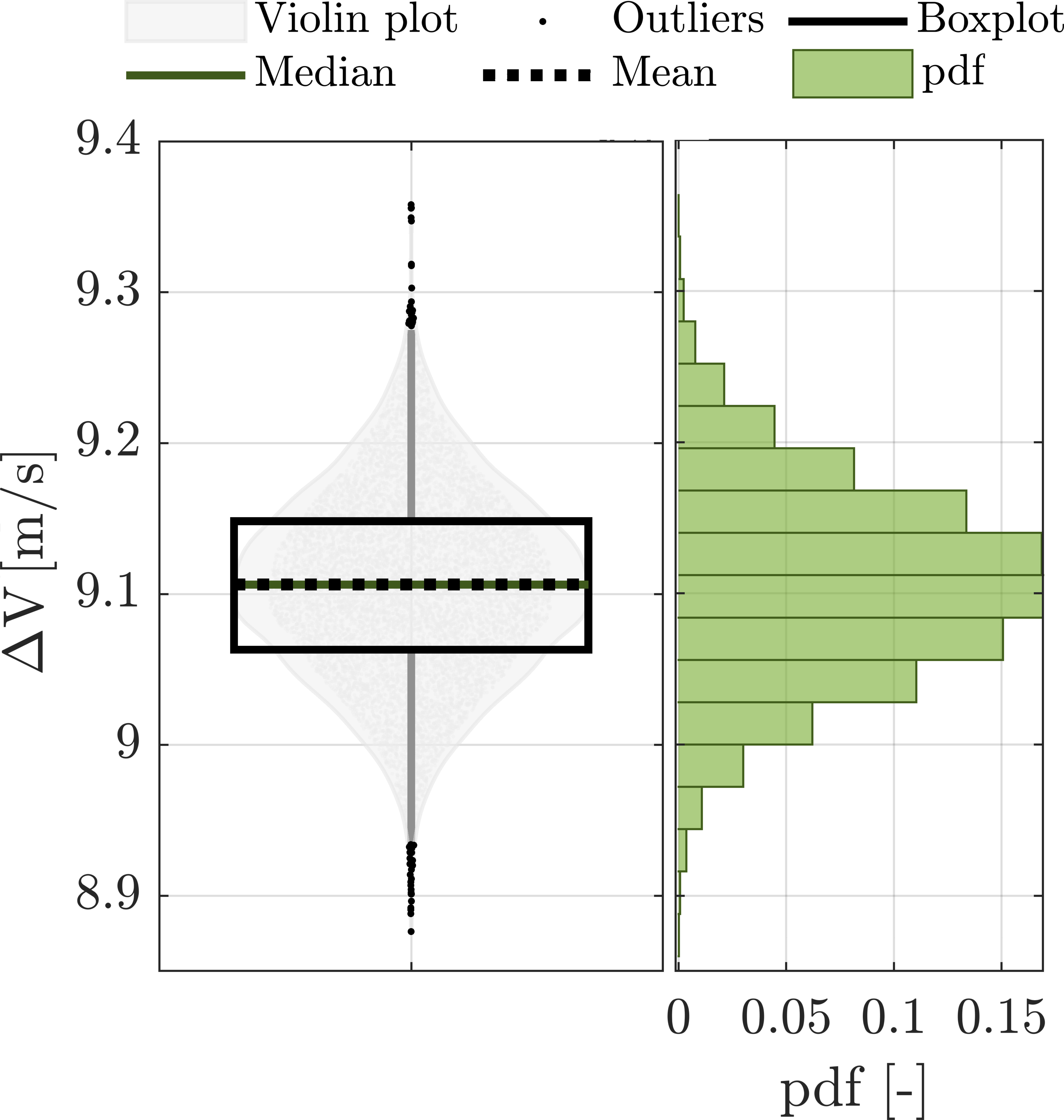}
  \end{center}
  \caption{A-posteriori probability density reconstruction by fitting a Pearson distribution in the first four moments. The maneuver magnitude can in first approximation be considered Gaussian. \label{fig::LTGTOGEO2}}
\end{figure}

\noindent \Cref{fig::LTGTOGEO2} shows the  reconstructed Pearson distribution's probability density, boxplot, and violin plot. This time \(s=0.019\) and \(k=2.994\) meaning the distribution can be, in first approximation, considered gaussian. The choice of using the \gls{cut}-4 algorithm thus allowed us to estimate different output statistics.  

\section{Conclusions}
\label{sec::conclusions}
This article introduced a novel method to detect and estimate a maneuver between two orbits exploiting convex optimization, with the only assumption that maneuvers were performed optimally. The minimum-fuel \gls{ocp} was transformed from a \gls{nlp} problem to a \gls{socp} problem to exploit the benefits of convex optimization: polynomial run-time complexity and convergence to the global optimum. The problem was discretized in time, \glspl{stm} were generated to linearly propagate the influence of maneuvers on the initial accurate trajectory, and a limit on the impulsive \(\Delta\)V was set to accommodate for different propulsion systems. In the end, it was possible to use successive convex optimization problems in case the optimum solution of the \gls{socp} was not the solution of the \gls{nlp} problem due to linearization. Results were shown for an out-of-plane maneuver in keplerian dynamics, an \gls{ewsk} maneuver in \gls{geo} and a \gls{lt} orbit raise from \gls{gto} to \gls{geo}, to highlight the versatility of the approach. \gls{cc} was fine to reconstruct the first two cases, but needed successive iterations to converge for the latter test, while \gls{coe} and \gls{mee} did not need any iteration, with the forward propagation of the maneuvering orbit always matching the \gls{socp} solution. Uncertainty in the initial and final states was dealt with through the Mahalanobis distance and \gls{cut}-4 approaches, which allowed for the detection of the maneuver and then estimation of mean, variance, skewness, and kurtosis parameters. By fitting a Pearson statistics through the four moments it was possible to show that the \(\Delta\)V statistics did not necessarily follow a Normal distribution. The tailorability of the appraoch was also shown by automatically restricting the search time window and using different dynamics. 

\section*{Aknowledgements}
The work presented is supported by AOARD under Grant FA2386-21-1-4115. The authors are grateful to Dr Jovan Skuljan of \gls{dta} for providing data and constructive discussions.

\bibliographystyle{unsrt} 
\bibliography{reference.bib}

\end{document}